\documentclass[12pt,fleqn,english,intoc,bibliography=totoc,BCOR10mm,captions=tableheading]{article}
\usepackage{lmodern}

\usepackage[T1]{fontenc}
\usepackage[utf8]{inputenc}
\usepackage{babel}
\usepackage{mathrsfs}
\usepackage{tipa}
\usepackage{amsmath}
\usepackage{amsthm}
\usepackage{amssymb}
\usepackage{stackrel}
\usepackage[a4paper]{geometry}
\usepackage{nomencl}
\makenomenclature
\usepackage[bookmarks=true,bookmarksnumbered=true,bookmarksopen=true,bookmarksopenlevel=1,
 breaklinks=true,pdfborder={0 0 1},backref=false,colorlinks=false]
 {hyperref}
\hypersetup{
 pdfauthor={Li, Wenwei},
 pdfborderstyle=,pdfpagelayout=OneColumn,pdfnewwindow=true,pdfstartview=XYZ,plainpages=false}

\makeatletter
\newcommand{\lyxaddress}[1]{
	\par {\raggedright #1
	\vspace{1.4em}
	\noindent\par}
}
\theoremstyle{plain}
\newtheorem{thm}{\protect\theoremname}
\theoremstyle{plain}
\newtheorem{lem}[thm]{\protect\lemmaname}

\usepackage{babel}
\usepackage{mathrsfs}
\usepackage{fancyhdr}
\usepackage{exscale}
\usepackage{relsize}
\usepackage{mathcomp}

\AtBeginDocument{
  
}

\providecommand{\lemmaname}{Lemma}
\providecommand{\theoremname}{Theorem}

\usepackage{ulem}
\makeatother

\providecommand{\lemmaname}{Lemma}
\providecommand{\theoremname}{Theorem}

\begin{document}
\title{\textsf{Some Relations on Paratopisms and an Intuitive Interpretation
on the Parastrophes of a Latin Square}\thanks{ Corresponding author. Wen-Wei Li (W.-W. Li), E-mail address: \emph{liwenwei@ustc.edu}.} }
\author{Wen-Wei Li $^{\mathrm{a,b,}}$, Jia-Bao Liu $^{\mathrm{c}}$, Xin
Hou $^{\mathrm{d}}$}
\maketitle

\lyxaddress{}

\lyxaddress{\begin{center}
{\footnotesize$^{\mathrm{a}}$ School of Mathematical Science, University
of Science and Technology of China, }\\
 {\footnotesize Hefei, 230026, P. R. China}\\
 {\footnotesize$^{\mathrm{b}}$ Research Center of Applied Mathematical
Models, Anhui International Studies University, }\\
 {\footnotesize Hefei, 231201, P. R. China}\\
 {\footnotesize$^{\mathrm{c}}$ School of Mathematics and Physics,
Anhui Jianzhu University, }\\
 {\footnotesize Hefei, 230601, P. R. China} \\
$^{\mathrm{d}}${\scriptsize{} }{\footnotesize College of Elementary
Education, Capital Normal University, }\\
{\footnotesize Beijing, 100048, P. R. China}
\par\end{center}}
\begin{abstract}
This paper  presents  some intuitive  interpretations  of the parastrophe
transformations of  an  arbitrary Latin square. With this trick, we
can generate the parastrophes of  the  arbitrary Latin square directly
from the original one without generating  an  orthogonal array. The
 relationships between  isotopisms and parastrophe transformations
in composition  are also  shown. It  solves  the problem that when
F1{*}I1=I2{*}F2 how can we obtain I2 and F2 from I1 and F1, where
I1 and I2 are isotopisms,  while F1 and F2 are parastrophe transformations,
 and ``{*}'' is the composition of transformations. These methods
could  significantly  simplify the computation for the issues related
to  the  main classes of Latin squares. This will improve the  computational
efficiency apparently in  some related problems.

\vspace{0.5cm}


\textbf{AMS2020 Subject Classification:} 05B15 
\end{abstract}
\newpage\tableofcontents{}

$\ $

\section{Introduction}

Latin squares play an important role in experimental design in combinatorics
and statistics. They are wildly used in  manufacturing industries,
agricultural productions, etc. Besides the application in  the  orthogonal
experiment in agriculture, Latin  squares can  also be used to minimize
 experimental  errors in the design of experiments in agronomic research
(refer  to  \cite{Nada2001ALSAR}).  A  Latin square is a multiplication
table of a quasi-group, which is an important structure in algebra.
Mutually orthogonal Latin squares are closely linked with finite projective
planes  \cite{Lam1991PP9B,Lam1991PP10}. Sets of orthogonal Latin
squares can be applied in error correcting codes in communication
(refer  to  \cite{Colbourn2004PAPCMOLS,Huczynska2006PC36OP}). It
is also used in Mathematical puzzles,  such as Sudoku  \cite{Ruiter2010-OJSPRT},
KenKen (refer  to  \cite{Shortz2009ANPCMS,Stephey2009IKKNS}), etc.

The name ``Latin square''  was  originated for the first time in
the 36 - officers problem introduced by Leonhard Euler (refer  to
 \cite{Euler1782RNEQM} or \cite{Euler2007IONTMS}). He used Latin
characters as symbols,  for which it was named as  ``Latin square''
 (instead of that, now we usually use the Hindu-Arabic numerals).
 Since  then, a lot of results have been obtained  in  this discipline.

One of the main  clues  in the development of this subject is the
enumeration problem, such as the number $L_{n}$ of reduced Latin
squares of order $n$ (a positive integer) or,  the number of equivalence
classes (such as isotopy classes or main classes) of Latin squares
of small orders. It is clear that the number $N_{n}$ of Latin squares
of order $n$ is $n!\cdot(n-1)!$ times of the number $L_{n}$ (of
the reduced Latin squares of order $n$). But there is no explicit
 relationship  (formula)  between  $L_{n}$ and the number $S_{n}^{(1)}$
of  the  isotopy classes of Latin squares of order $n$, which  can
 be used to compute $S_{n}^{(1)}$ from $L_{n}$ or $n$ in  practice.
 The number $S_{n}^{(2)}$ of the main classes of Latin squares of
order $n$  cannot  be calculated directly from $L_{n}$ or $n$,
either.

Till today, there is no practical  formula for computing  $L_{n}$,
 which  can  easily obtain $L_{n}$.  Shao and Wei   derived  a simple
and explicit formula (in form)  in 1992 (refer to \cite{ShaoJY1992AFNLS})
for computing for $L_{n}$ as  $L_{n}=n!\underset{A\in B_{n}}{\sum}(-1)^{\sigma_{0}(A)}\dbinom{\mathrm{Per}A}{n}$,
where $B_{n}$ is the set of all the 0-1 square matrices of order
$n$, $\sigma_{0}(A)$ is the number of ``0'' appeared in matrix
$A$,  and  ``$\mathrm{Per}$'' is the permanent operator. However,
 this formula is still not so  efficient  in  practice. There is no
practical asymptotic formula for  computing  $L_{n}$, either (refer
 to  \cite{McKay2005ONLS}). The difference  between  the most accurate
upper bounds and lower bounds of $L_{n}$ is huge, $\dfrac{\left(n!\right)^{2n}}{n^{n^{2}}}\leqslant L_{n}\leqslant\stackrel[k=1]{n}{\prod}\left(k!\right)^{n/k}$
( as  mentioned in \cite{Lint1992ACC}, pp.161-162), which made it
impossible to estimate the value of $L_{n}$ by  using  this formula.
The upper bound was inferred from the  permanent conjecture of van
der Waerden  by Richard M. Wilson in 1974  \cite{Wilson1974NSTS}.
Later  on, in  around 1980, the van der Waerden conjecture was solved
independently by G. P. Egorichev \cite{Egorichev1981PVDWCRPDSM-Eng,Egorichev1981PVDWCP-Eng}
 and D. I. Falikman (refer  to  \cite{Falikman1981PVDWCRPDSM-Rus,Falikman1981PVDWCRPDSM-Eng}),
 almost at the same time.

James R. Nechvatal ( as  mentioned in \cite{Godsil1990AEGLR})  presented
in 1981  a general asymptotic formula for generalized Latin Rectangles,
and Ira M. Gessel \cite{Gessel1987CLR}   presented in 1987  a general
asymptotic formula for Latin Rectangles. Although  they could  derive
theoretical  asymptotic formulae of $L_{n}$,  none was found to be
 suitable for asymptotic analysis. Chris D. Godsil and Brendan D.
McKay  derived  a better asymptotic formula in 1990  \cite{Godsil1990AEGLR}.

When the order $n$ is less than 4, the number $L_{n}$ of reduced
Latin squares of order $n$ is obvious, i.e., $L_{1}=L_{2}=L_{3}=1$.
 For  $n=4$ or 5, Euler found  in 1782  that $L_{4}=4$ and $L_{5}=56$
(refer  to  \cite{Euler1782RNEQM}), together with the values of $L_{1}$,
$L_{2}$ and $L_{3}$. Cayley  also found  these results (up to 5)
in 1890 (refer  to  \cite{Cayley1890OLS}). M. Frolov found  $L_{6}=9,408$
 in 1890 ( as  mentioned in \cite{Brendan2007SmallLS}). Later on,
Tarry  re-found it in 1901  ( as  mentioned in \cite{Brendan2007SmallLS}).
The number $S_{6}^{(1)}=22$ was  obtained first  by E. Schönhardt
 \cite{Schoenhardt1930ULQU} in 1930. R. A. Fisher and F. Yates  \cite{Fisher1934LSOd6}
 also found  $S_{6}^{(1)}$ independently in 1934,  as well as  the
values of $S_{n}^{(1)}$  for  $n\leqslant5$. In 1966, D. A. Preece
 \cite{Preece1966CYR} found that there are 564 isotopy classes of
Latin squares of order 7. M. B. Wells  \cite{Wells1967NLSOd8} acquired
$L_{8}=535,281,401,856$ in 1967. In 1990, G. Kolesova et al.  \cite{Lam1990NumLS8}
 found  $S_{8}^{(1)}=$ 1,676,267, $S_{8}^{(2)}$ = 283,657, which
 confirmed  the  Wells' result. The value of $L_{9}=$ 377,597,570,964,258,816
was found by S. E. Bammel and J. Rothstein  \cite{Bammel1975NLSOd9}
in 1975. B. D. McKay and E. Rogoyski  \cite{McKayl1995LSOd10} found
 in 1995  that $L_{10}=$ 7,580,721,483,160,132,811,489,280,  $S_{10}^{(1)}$
= 208,904,371,354,363,006 and $S_{10}^{(2)}$ = 34,817,397,894,749,939.
B. D. McKay and I. M. Wanless  \cite{McKay2005ONLS} found  in 2005
 the  values  of $L_{11}$ = 5,363,937,773,277,371,298,119,673,540,771,840,
 $S_{11}^{(1)}$ =12,216,177,315,369,229,261,482,540 and $S_{11}^{(2)}$
= 2,036,029,552,582,883,134,196,099.

$\ $

Before the invention of computers,  mathematicians used  do the enumeration
work by hand with some theoretical tools, which would probably involve
some errors. It was also very difficult to verify a known result.
Hence,  some experts obtained false values even if the  correct  one
had been found. In 1915, Percy A. MacMahon  obtained  an incorrect
value of $L_{5}$ in a different way from other experts (refer  to
 \cite{MacMahon1915CombAnl}). In 1930, S. M. Jacob  (\cite{Jacob1930ELRD3})
obtained a wrong value of $L_{6}$  even  after Frolov and Tarry had
already found the corrected one. The value of $L_{7}$  calculated
 by Frolov is  incorrect.

It was  mentioned  in \cite{Norton1939LSOd7} that Clausen, an assistant
of a German astronomer Schumacher, found 17 ``basic forms''\footnote{$\ $ Translated from the German word ``Grundformen''. According
to the context, it probably means the isotopy classes. } of Latin squares of order 6. This  observation  was described in
a letter from Schumacher to Gauss dated August 10, 1842. This letter
was quoted by Gunther in 1876 and by Ahrens in 1901. Tarry also found
17 isotopy classes of Latin squares of order 6. (E.  Schönhardt obtained
 the correct values of $L_{n}$, $S_{n}^{(1)}$ and $S_{n}^{(2)}$
 for  $n\leqslant6$ in 1930.)

In 1939, H. W. Norton \cite{Norton1939LSOd7}  obtained  some  wrong
values of $S_{7}^{(1)}$ and $S_{7}^{(2)}$. After Preece  found  the
correct value of $S_{7}^{(1)}$ in 1966,  Brown  \cite{Brown1968ELSAOd8}
 reported  another incorrect value  of $S_{7}^{(1)}=563$. This  result
was widely quoted as the accepted value for several decades  \cite{Colbourn1996CRCHBCD,Denes1974Latin}.

In recent decades, with the application of computers, the efficiency
in  enumerating  equivalence classes of Latin squares has been improved
greatly. However, it is still very difficult to avoid errors because
of the complexity  of the involved  huge amount of computation.

In the case of order 8, J. W. Brown  \cite{Brown1968ELSAOd8} also
 presented  a wrong value of $S_{8}^{(1)}$ in 1968, and Arlazarov
et al. provided a false value of $S_{8}^{(2)}$ in 1978 ( as  mentioned
in \cite{Brendan2007SmallLS,Lam1990NumLS8}).

In some cases, the number of some types of equivalence classes of
Latin squares of a certain order would likely to be believed correct
after at least two times of independent computation  (refer to \cite{Brendan2007SmallLS}). 

$\ $

Technically, according to some conclusions in  the  group theory,
especially the enumeration method of orbits  when a group acts on
a set,  we can know the number of objects in every equivalence class
by generating representatives of all the equivalence classes of a
certain order and then enumerating the members in each invariant group
of a representative.   Accordingly,  the total number of Latin squares
of a certain order  can  be obtained. But the process contains  many
 equivalence classes, which costs  a lot of  time.\footnote{$\ $ The number of isotopy classes or main classes of Latin squares
of order $n$ is much less than the number of reduced Latin squares
of order $n$. In general, we  cannot  afford the time to visit all
the reduced Latin squares of order  $n>7$  even in some super computers.}

 Logically, the structure of the relations of the Latin squares in
the same main class is a graph. The  powerful graph isomorphism program
``Nauty'' is very useful  in computing  the invariant group of  a
 Latin square within  the  main class transformations.

In  practice,  the process of computing these objects (to find the
invariant group of a Latin square) is  similar to visiting a tree
with some of its branches being isomorphic,  where only one of the
isomorphic branches   should be visited  for improving  the efficiency.

When considering the topics related to the main classes of Latin squares,
we usually generate the parastrophes (or conjugates) of an arbitrary
Latin square. The routine procedure is to generate the orthogonal
array of the Latin square first, then  to  permute the rows (or columns
in some literature) of the array,  and  finally to  turn the new orthogonal
array into   a new Latin square, which seems not  very  convenient.
Among the parastrophes of a Latin square $Y$, the authors have not
found the detailed descriptions  of  the parastrophes other than the
orthogonal array, except two simple cases  of  $Y$ itself and  its
 transpose $Y^{\mathrm{T}}$. In Sec. \ref{subsec:Adjugates-LS},
an  intuitive explanation of the parastrophes of an arbitrary Latin
square will be described.  Using  this method, other 4 types  of  parastrophes
of  an  arbitrary Latin square  can  be generated by  the  transpose,
 and/or  by  replacing the rows/columns  with  their inverse.

Sometimes we need to consider the composition of an parastrophe transformation
and an isotopic transformation, especially when generating the invariant
group of a Latin square in  the  main class transformation. It is
necessary to exchange the priority order of the these two types of
transformations for simplification. But in general cases, an parastrophe
transformation and an isotopic transformation do not commute. So,
 we  need  to find the relations  between the  isotopic transformation
and the parastrophe transformation in  the  composition when their
positions are interchanged. In other words, for any parastrophe transformation
$\mathcal{F}_{1}$ and isotopic transformation $\mathcal{I}_{1}$,
how can we find the parastrophe transformation $\mathcal{F}_{2}$
and isotopic transformation $\mathcal{I}_{2}$,  s.t., $\mathcal{F}_{1}\circ\mathcal{I}_{1}$
= $\mathcal{I}_{2}\circ\mathcal{F}_{2}$ ? The answer will be presented
in Sec. \ref{subsec:RLT-Prtp}.

\section{Preliminaries}

 Some notions used in this paper are defined here for avoiding any
ambiguity. 

Suppose $n$  is a positive integer and it is greater than 1.

A \emph{permutation}\index{permutation} \label{Def:Permutation}
is the reordering of the sequence  of  1, 2, 3, $\cdots$, $n$. An
element $\alpha$ = $\left(\begin{array}{cccc}
1 & 2 & \cdots & n\\
a_{1} & a_{2} & \cdots & a_{n}
\end{array}\right)$ in the symmetry group $\mathrm{S}_{n}$ is also called a \emph{permutation}.
\label{def:permutation}

For convenience, in this paper, these two concepts will not be  differentiated
 rigorously. They might even be  used interchangeably. When referring
 to  ``a permutation $\alpha$'' here, sometimes it is  a bijection
of the set \{$\,$1, 2, 3, $\cdots$, $n$$\,$\} to itself,  while
 sometimes it  stands for the sequence $\left[\alpha(1),\,\alpha(2),\,\cdots,\,\alpha(n)\right]$.
For a sequence $\left[b_{1},b_{2},\cdots,b_{n}\right]$,  which is
a rearrangement of {[}1, 2, $\cdots$, $n${]}, in some occasions
it may also stand for a transformation $\beta\in\mathrm{S}_{n}$,
 so  that $\beta(i)=b_{i}$ ($i$ = 1, 2, $\cdots$, $n$).  Its  actual
meaning can be inferred from the  associated  contexts.\label{CVT:Permutation-2}
In  a  computer, these two  types  of objects are stored almost in
the same way.

We are compelled to accept this ambiguity. Otherwise, it will cost
too much energy to avoid this ambiguity  as  we have to use much more
words to describe a simple operation and much more symbols to show
a concise expression. An example will be shown after Lemma 1.

Let $\alpha$ = $\left(\begin{array}{cccc}
1 & 2 & \cdots & n\\
a_{1} & a_{2} & \cdots & a_{n}
\end{array}\right)$ $\in\mathrm{S}_{n}$,  where $\left[a_{1},a_{2},\cdots,a_{n}\right]$
is  the \emph{one-row form} of the permutation $\alpha$, and  $\left(\begin{array}{cccc}
1 & 2 & \cdots & n\\
a_{1} & a_{2} & \cdots & a_{n}
\end{array}\right)$ is  the \emph{two-row form} of  $\alpha$.

A matrix with  its  every row and every column being  a permutation
 of 1, 2, $\cdots$, $n$,\footnote{$\ $  For convenience, we will usually assume that the $n$ elements
of a Latin square are 1, 2, 3, $\cdots$, $n$.  But in a lot of books
and articles, the $n$ elements  of  a Latin  square  are denoted
by 0, 1, 2, 3, $\cdots$, $n-1$.} is called a \emph{Latin square}\index{Latin square} of order $n$.
Latin squares with both the first row and first column being in natural
order are said to be \emph{reduced}\index{reduced (Latin square)}
or in \emph{standard form}\index{standard form} (refer  to  \cite{Denes1974Latin},
 page 128). For example, the  following matrix  is a reduced Latin
square of order 5 :  
\[
\left[\begin{array}{ccccc}
1 & 2 & 3 & 4 & 5\\
2 & 3 & 5 & 1 & 4\\
3 & 5 & 4 & 2 & 1\\
4 & 1 & 2 & 5 & 3\\
5 & 4 & 1 & 3 & 2
\end{array}\right].
\]

For convenience, when referring  to  the \emph{inverse} of a row (or
a column) of a Latin square, we mean the one-row form of the inverse
of the permutation presented by the row (or the column), not the sequence
in  a  reverse order. For instance, the inverse of the third row $[3\ 5\ 4\ 2\ 1]$
of the  above Latin square  is believed to be $[5\ 4\ 1\ 3\ 2]$,
not  its reverse $[1\ 2\ 4\ 5\ 3]$ as  $\left(\begin{array}{ccccc}
1 & 2 & 3 & 4 & 5\\
3 & 5 & 4 & 2 & 1
\end{array}\right)^{-1}=$ $\left(\begin{array}{ccccc}
1 & 2 & 3 & 4 & 5\\
5 & 4 & 1 & 3 & 2
\end{array}\right)$.\footnote{$\ $ It is not difficult to find out that, in computer programs,
if we store a permutation $\left[a_{1},a_{2},\cdots,a_{n}\right]$
in an array ``A{[}n{]}'', the inverse ``B{[}n{]}'' of ``A{[}n{]}''
 can  be generated  by evaluating \emph{n} times,  ``B{[}A{[}j{]}{]}=j'',
j = 1, 2, $\cdots$, $n$. The traditional method used by some programmers
 interchanges  the rows of the permutation $\left(\begin{array}{ccccc}
1 & 2 & 3 & 4 & 5\\
3 & 5 & 4 & 2 & 1
\end{array}\right)$,  and  then  sorts  the columns of the new permutation $\left(\begin{array}{ccccc}
3 & 5 & 4 & 2 & 1\\
1 & 2 & 3 & 4 & 5
\end{array}\right)$  to make the first row in  natural order, so as to obtain the inverse
$\left(\begin{array}{ccccc}
1 & 2 & 3 & 4 & 5\\
5 & 4 & 1 & 3 & 2
\end{array}\right)$. The number of operations  required in the  evaluations and comparisons
will be much  higher  than  that required in  the previous method. }

Let $\forall\alpha,\beta,\gamma\in\mathscr{\mathrm{S}}_{n}$,  $Y=\left(y_{ij}\right)_{n\times n}$
be an arbitrary Latin square with  its  elements belonging to the
set  of  $\left\{ \,1,\,2,\,\mbox{\ensuremath{\cdots},\,}n\,\right\} $. 

$\mathscr{R}_{\alpha}$ \label{Sym:mathscr_R_beta} \nomenclature[R_alpha]{$\mathscr{R}_{\beta}$}{the transformation that {permutes} the rows of a Latin square according to a permutation $\beta$. \pageref{Sym:mathscr_R_beta}}
 denotes  the transformation of rows  in  a Latin square (or a matrix)
corresponding to the permutation $\alpha$, i.e., the $\alpha(i)$'s
row of the Latin square $\mathscr{R}_{\alpha}\left(Y\right)$ is the
$i$'s row of the Latin square $Y$, in other words, the $i$'s row
of the Latin square $\mathscr{R}_{\alpha}\left(Y\right)$ is the $\alpha^{-1}(i)$'s
row of the Latin square $Y$. $\mathscr{C}_{\beta}$ \label{Sym:mathscr_C_alpha}
\nomenclature[C_\alpha]{$ \mathscr{C}_{\alpha}$}{the transformation that {permutes} the columns of a Latin square according to a permutation $\alpha$ of the same order.  \pageref{Sym:mathscr_C_alpha}}
 denotes  the transformation  that  moves the $i$'th column ($i$
= 1, 2, $\cdots$, $n$) of a Latin square (or a general matrix) to
the position of the $\beta(i)$'th column. $\mathscr{L}_{\gamma}$
\label{Sym:mathscr_L_gamma} \nomenclature[L_gamma]{$\mathscr{L}_{\gamma}$}{the transformation that substitutes the elements of a Latin square according to a permutation $\gamma$. \pageref{Sym:mathscr_L_gamma}}
 denotes  the transformation  that  relabels the entries, i.e., $\mathscr{L}_{\gamma}$
 substitutes  all the elements $i$ in a Latin square by $\gamma(i)$
($i$ = 1, 2, $\cdots$, $n$). For example, let $\alpha$=$\left(\begin{array}{cccc}
1 & 2 & 3 & 4\\
2 & 4 & 1 & 3
\end{array}\right)$, $\beta$=$\left(\begin{array}{cccc}
1 & 2 & 3 & 4\\
2 & 3 & 4 & 1
\end{array}\right)$, $\gamma$=$\left(\begin{array}{cccc}
1 & 2 & 3 & 4\\
3 & 1 & 4 & 2
\end{array}\right)$,  and  $Y_{1}$=$\left[\begin{array}{cccc}
1 & 2 & 3 & 4\\
2 & 3 & 4 & 1\\
3 & 4 & 1 & 2\\
4 & 1 & 2 & 3
\end{array}\right]$. Then  $\mathscr{R}_{\alpha}\left(Y_{1}\right)$ = $\left[\begin{array}{cccc}
3 & 4 & 1 & 2\\
1 & 2 & 3 & 4\\
4 & 1 & 2 & 3\\
2 & 3 & 4 & 1
\end{array}\right]$, $\mathscr{C}_{\beta}\left(Y_{1}\right)$ = $\left[\begin{array}{cccc}
4 & 1 & 2 & 3\\
1 & 2 & 3 & 4\\
2 & 3 & 4 & 1\\
3 & 4 & 1 & 2
\end{array}\right]$,  and  $\mathscr{L}_{\gamma}\left(Y_{1}\right)$ = $\left[\begin{array}{cccc}
3 & 1 & 4 & 2\\
1 & 4 & 2 & 3\\
4 & 2 & 3 & 1\\
2 & 3 & 1 & 4
\end{array}\right]$.

\label{CVT:Permutation-3}  Now onward, let $Y_{i}$ = $\left[\begin{array}{cccc}
y_{i1}, & y_{i2}, & \cdots, & y_{in}\end{array}\right]$ be the $i$'th row of a Latin square $Y=\left(y_{ij}\right)_{n\times n}$,
and  $Z_{i}$ = $\left[\begin{array}{c}
y_{1i}\\
y_{2i}\\
\vdots\\
y_{ni}
\end{array}\right]$ be the $i$'th column of $Y$ ($i$=1, 2, $\cdots$, $n$).  Define
two transformations $\eta_{i}$ and $\zeta_{i}$, where $\eta_{i}$
: \{1,\emph{ }2, $\cdots$\emph{,} $n$\} $\rightarrow$ \{1, 2, $\cdots$,
$n$\}, $j$ $\longmapsto$ $y_{ij}$\emph{ }($j$=1, 2, $\cdots$,
$n$);  $\zeta_{i}$ : \{1,\emph{ }2, $\cdots$\emph{,} $n$\} $\rightarrow$
\{1, 2, $\cdots$, $n$\}, $k$ $\longmapsto$ $y_{ki}$\emph{ }($k$=1,
2, $\cdots$, $n$). By the assumption mentioned after the definition
of permutations on page \pageref{Def:Permutation}, we do not distinguish
the permutation transformation $\zeta$ : \{1,\emph{ }2, $\cdots$\emph{,}
$n$\} $\rightarrow$ \{1, 2, $\cdots$, $n$\}, $k$ $\longmapsto$
$\zeta(k)$\emph{ }($k$=1, 2, $\cdots$, $n$) and the column sequence
{[}$\zeta(1)$,\emph{ }$\zeta(2)$, $\cdots$\emph{,} $\zeta(n)${]}$^{\mathrm{T}}$.
Here,  the super-script ``T'' means transpose. So, $Y_{i}$ = $\left[\begin{array}{cccc}
y_{i1}, & y_{i2}, & \cdots, & y_{in}\end{array}\right]$ and $\eta_{i}$ = $\left(\begin{array}{cccc}
1 & 2 & \cdots & n\\
y_{i1} & y_{i2} & \cdots & y_{in}
\end{array}\right)$ will be considered as the same object; $Z_{i}$ = $\left[\begin{array}{c}
y_{1i}\\
y_{2i}\\
\vdots\\
y_{ni}
\end{array}\right]$ and $\zeta_{i}$ = $\left(\begin{array}{cccc}
1 & 2 & \cdots & n\\
y_{1i} & y_{2i} & \cdots & y_{ni}
\end{array}\right)$ will not be differentiated.  $Y_{i}^{-1}$ and $Z_{i}^{-1}$ are
the  inverses  of $Y_{i}$ and $Z_{i}$, respectively, i.e., $Y_{i}^{-1}$
and $Z_{i}^{-1}$ are regarded as the one-row form of the transformations
$\eta_{i}^{-1}$ = $\left(\begin{array}{cccc}
y_{i1} & y_{i2} & \cdots & y_{in}\\
1 & 2 & \cdots & n
\end{array}\right)$ and $\zeta_{i}^{-1}$ = $\left(\begin{array}{cccc}
y_{1i} & y_{2i} & \cdots & y_{ni}\\
1 & 2 & \cdots & n
\end{array}\right)$, respectively. Besides, $Y_{i}^{-1}$ is a row and $Z_{i}^{-1}$
is a column. \footnote{$\ $ Actually, $Z_{i}^{-1}$ = $\left[\begin{array}{c}
\zeta_{i}^{-1}(1)\\
\zeta_{i}^{-1}(2)\\
\vdots\\
\zeta_{i}^{-1}(n)
\end{array}\right]$ is the transpose  of  the sequence {[}$\zeta_{i}^{-1}(1)$,\emph{
}$\zeta_{i}^{-1}(2)$, $\cdots$\emph{,} $\zeta_{i}^{-1}(n)${]},
where $\zeta_{i}^{-1}$ is the inverse of the permutation $\zeta_{i}$
= $\left(\begin{array}{cccc}
1 & 2 & \cdots & n\\
y_{1i} & y_{2i} & \cdots & y_{ni}
\end{array}\right)$. $Y_{i}^{-1}$ is obtained in the same way, except that $Y_{i}^{-1}$
is a row sequence.} Whether a permutation symbol $\alpha$ stands for a row sequence
or  a  column sequence will be inferred from the  corresponding  contexts.

For $\forall\alpha,\beta,\gamma\in\mathscr{\mathrm{S}}_{n}$, we know
the definition of the composition $\alpha\cdot\beta$ as a permutation
transformation,  where  $\left(\alpha\cdot\beta\right)(i)$ is defined
 as  $\alpha\bigl(\beta(i)\bigr)$, $i$=1, 2, $\cdots$, $n$. We
define $\gamma Y_{i}$ as  the one-row form of the composition $\gamma\cdot\eta_{i}$
of permutations $\gamma$ and $\eta_{i}$, where $\eta_{i}$ is the
permutation  $\left(\begin{array}{cccc}
1 & 2 & \cdots & n\\
y_{i1} & y_{i2} & \cdots & y_{in}
\end{array}\right)$,  i.e, usually, $\gamma Y_{i}$  stands  for the sequence  of  $\left[\begin{array}{cccc}
\gamma\left(y_{i1}\right), & \gamma\left(y_{i2}\right), & \cdots, & \gamma\left(y_{in}\right)\end{array}\right]$, except that in a few occasions $\gamma Y_{i}$ is the permutation
transformation $\gamma\cdot\eta_{i}$ = $\left(\begin{array}{cccc}
1 & 2 & \cdots & n\\
\ \gamma\left(y_{i1}\right)\  & \ \gamma\left(y_{i1}\right)\  & \cdots & \ \gamma\left(y_{i1}\right)\ 
\end{array}\right)$, according to the  corresponding contexts.  So are $Z_{i}\alpha^{-1}$
and $\gamma Z_{i}$, except that $Z_{i}$ is a column, hence both
$Z_{i}\alpha^{-1}$ and $\gamma Z_{i}$ are columns. 

It is not difficult to verify  the following Lemma:  
\begin{lem}
\label{lem:Act-Istp} $\ $ Let $\forall\alpha,\beta,\gamma\in\mathscr{\mathrm{S}}_{n}$,
 $Y=\left(y_{ij}\right)_{n\times n}$ = $\left(\begin{array}{c}
Y_{1}\\
\vdots\\
Y_{n}
\end{array}\right)$ = $\left(\begin{array}{ccc}
Z_{1}, & \cdots, & Z_{n}\end{array}\right).$ be a Latin square of order $n$, where $Y_{i}$ is the $i$'th row
of $Y$ and $Z_{i}$ is the $i$'th column of $Y$ \emph{(}$i$=1,
2, $\cdots$, $n$\emph{)}, then
\begin{align*}
\mathscr{R}_{\alpha}(Y) & =\left(\begin{array}{c}
Y_{\alpha^{-1}(1)}\\
\vdots\\
Y_{\alpha^{-1}(n)}
\end{array}\right)=\left(\begin{array}{ccc}
Z_{1}\alpha^{-1}, & \cdots, & Z_{n}\alpha^{-1}\end{array}\right),\\
\mathscr{C}_{\beta}(Y) & =\left(\begin{array}{c}
Y_{1}\beta^{-1}\\
\vdots\\
Y_{n}\beta^{-1}
\end{array}\right)=\left(\begin{array}{ccc}
Z_{\beta^{-1}(1)}, & \cdots, & Z_{\beta^{-1}(n)}\end{array}\right),\\
\mathscr{L}_{\gamma}(Y) & =\left(\begin{array}{c}
\gamma Y_{1}\\
\vdots\\
\gamma Y_{n}
\end{array}\right)=\left(\begin{array}{ccc}
\gamma Z_{1}, & \cdots, & \gamma Z_{n}\end{array}\right).
\end{align*}
\end{lem}

If we distinguish the two  types  of permutations strictly, it will
make the above equation much more complex. For example, in order to
show that the transformation  of  $\left(\begin{array}{cccc}
1 & 2 & \cdots & n\\
a_{1} & a_{2} & \cdots & a_{n}
\end{array}\right)$ is derived from a sequence  of  $R=\left[a_{1},a_{2},\cdots,a_{n}\right]$,
(which  is a reordering of {[}1, 2, $\cdots$, $n${]}), we will use
a symbol,  such as $\mathscr{S}(R)$,  to denote the transformation
 of  $\left(\begin{array}{cccc}
1 & 2 & \cdots & n\\
a_{1} & a_{2} & \cdots & a_{n}
\end{array}\right)$. For a transformation $\beta$ = $\left(\begin{array}{cccc}
1 & 2 & \cdots & n\\
b_{1} & b_{2} & \cdots & b_{n}
\end{array}\right)$, and sequence $\left[b_{1},b_{2},\cdots,b_{n}\right]$,  we should
denote the sequence $\left[b_{1},b_{2},\cdots,b_{n}\right]$ by $\mathscr{T}(\beta)$
in order to show the relationship between $\beta$ and the corresponding
sequence. Therefore,  Lemma 1 should be stated as  follows: 
\begin{align*}
\mathscr{R}_{\alpha}(Y) & =\left(\begin{array}{c}
Y_{\alpha^{-1}(1)}\\
\vdots\\
Y_{\alpha^{-1}(n)}
\end{array}\right)=\left(\begin{array}{ccc}
\left(\mathscr{T}\left(\mathscr{S}\left(Z_{1}^{\mathrm{T}}\right)\cdot\alpha^{-1}\right)\right)^{\mathrm{T}}, & \cdots, & \left(\mathscr{T}\left(\mathscr{S}\left(Z_{n}^{\mathrm{T}}\right)\cdot\alpha^{-1}\right)\right)^{\mathrm{T}}\end{array}\right),\\
\mathscr{C}_{\beta}(Y) & =\left(\begin{array}{c}
\mathscr{T}\left(\mathscr{S}\left(Y_{1}\right)\cdot\beta^{-1}\right)\\
\vdots\\
\mathscr{T}\left(\mathscr{S}\left(Y_{n}\right)\cdot\beta^{-1}\right)
\end{array}\right)=\left(\begin{array}{ccc}
Z_{\beta^{-1}(1)}, & \cdots, & Z_{\beta^{-1}(n)}\end{array}\right),\\
\mathscr{L}_{\gamma}(Y) & =\left(\begin{array}{c}
\mathscr{T}\left(\mathscr{\gamma\cdot S}\left(Y_{1}\right)\right)\\
\vdots\\
\mathscr{T}\left(\mathscr{\gamma\cdot S}\left(Y_{n}\right)\right)
\end{array}\right)=\left(\begin{array}{ccc}
\left(\mathscr{T}\left(\gamma\cdot\mathscr{S}\left(Z_{1}^{\mathrm{T}}\right)\right)\right)^{\mathrm{T}}, & \cdots, & \left(\mathscr{T}\left(\gamma\cdot\mathscr{S}\left(Z_{n}^{\mathrm{T}}\right)\right)\right)^{\mathrm{T}}\end{array}\right).
\end{align*}

Here,  the two transformation symbols  of  $\mathscr{S}$ and $\mathscr{T}$
will make the above equations not as concise as the previous ones,
and they will distract our attention. Also, it will cause  a big trouble
 in demonstrating  some other propositions which are more complicated
than Lemma 1. We will pay  too much attention to strictly distinguish
these two types of permutations. The price is too much. Therefore,
 we will ignore it in this paper.

Later, the two symbols  of  $\mathscr{S}$ and $\mathscr{T}$ will
stand for some other injective maps.

For $\forall\alpha,\beta,\gamma\in\mathscr{\mathrm{S}}_{n}$ and any
Latin square $Y$, $\mathscr{R}_{\alpha}\circ\mathscr{C}_{\beta}\circ\mathscr{L}_{\gamma}$
will be called an \emph{isotopism},  and the Latin square $H=\left(\mathscr{R}_{\alpha}\circ\mathscr{C}_{\beta}\circ\mathscr{L}_{\gamma}\right)(Y)$
will be called \emph{isotopic} to $Y$.

Let $\mathscr{I}_{n}$ be the set of all the isotopy transformations
of Latin squares of order $n.$

For a Latin square $Y=\left(y_{ij}\right)_{n\times n}$ with $y_{ij}\in\left\{ \,1,\,2,\,\mbox{\ensuremath{\cdots},\,}n\,\right\} $,
with regard to the set 
\[
\mathrm{T}=\left\{ \left.(i,j,y_{ij})\,\right|\,1\leqslant i,j\leqslant n\,\right\} ,
\]
 we have 
\[
\left\{ \left.(i,j)\,\right|\,1\leqslant i,j\leqslant n\,\right\} =\left\{ \left.(j,y_{ij})\,\right|\,1\leqslant i,j\leqslant n\,\right\} =\left\{ \left.(i,y_{ij})\,\right|\,1\leqslant i,j\leqslant n\,\right\} .
\]
 Hence,  each pair of triplets $(i,j,y_{ij})$ and $(r,t,y_{rt})$
in $\mathrm{T}$ will share at most one identical entry in the same
position. The set $\mathrm{T}$ is also called the \emph{orthogonal
array representation} \index{orthogonal array representation} of
the Latin square $Y$.

 Now onward, each  triplet  $(i,j,y_{ij})$ in the orthogonal array
set $\mathrm{T}$ of Latin square $Y=\left(y_{ij}\right)_{n\times n}$
will be written in the form of column vector $\left[\begin{array}{c}
i\\
j\\
y_{ij}
\end{array}\right]$ so as to save some space ( in many papers, a row vector is used to
represent the triplet ). The orthogonal array set $\mathrm{T}$ of
Latin square $Y$ can be  written in a matrix  
\begin{equation}
V=\left[\begin{array}{llll|llll|lll|llll}
1 & 1 & \cdots & 1 & \,2 & 2 & \cdots & 2 & \,\bullet & \bullet & \bullet\, & \,n & n & \cdots & n\\
1 & 2 & \cdots & n & \,1 & 2 & \cdots & n & \,\bullet & \bullet & \bullet & \,1 & 2 & \cdots & n\\
y_{11} & y_{12} & \cdots & y_{1n} & \,y_{21} & y_{22} & \cdots & y_{2n} & \,\bullet & \bullet & \bullet & \,y_{n1} & y_{n2} & \cdots & y_{nn}
\end{array}\right]\label{eq:OrthogonalArray}
\end{equation}
of size $3\times n^{2}$, with every column being a  triplet consisting
 of the indices of a position in a Latin square and the element in
that position.  Matrix  $V$ will also be called  as  the \emph{orthogonal
array }(matrix)\index{orthogonal array (matrix)}.   From now on,
when referring  to  the orthogonal array of a Latin square, it will
always be the matrix with every column being a  triplet  related to
a position of the Latin square, unless otherwise specified.

The definition of ``orthogonal array” in general may be found in
reference \cite{Denes1974Latin} (page 190).\emph{ }In some references,
such as \cite{Cameron2004LatinSquares}, the orthogonal array is defined
as an $n^{2}\times3$ array.

Every row of the orthogonal array of a Latin square consists of the
elements  of  1, 2, $\cdots$, $n$, and every element appears  exactly
 $n$ times in every row.

For example, the orthogonal arrays of the Latin squares $A_{1}=\left[\begin{array}{ccccc}
1 & 2 & 3 & 4 & 5\\
2 & 3 & 4 & 5 & 1\\
3 & 4 & 5 & 1 & 2\\
4 & 5 & 1 & 2 & 3\\
5 & 1 & 2 & 3 & 4
\end{array}\right]$ are \label{Eg:Orth-Arr} 
\[
V_{1}=\left[\begin{array}{ccccc|ccccc|ccccc|ccccc|ccccc}
1 & 1 & 1 & 1 & 1 & 2 & 2 & 2 & 2 & 2 & 3 & 3 & 3 & 3 & 3 & 4 & 4 & 4 & 4 & 4 & 5 & 5 & 5 & 5 & 5\\
1 & 2 & 3 & 4 & 5 & 1 & 2 & 3 & 4 & 5 & 1 & 2 & 3 & 4 & 5 & 1 & 2 & 3 & 4 & 5 & 1 & 2 & 3 & 4 & 5\\
1 & 2 & 3 & 4 & 5 & 2 & 3 & 4 & 5 & 1 & 3 & 4 & 5 & 1 & 2 & 4 & 5 & 1 & 2 & 3 & 5 & 1 & 2 & 3 & 4
\end{array}\right].
\]

Obviously,  the reordering of  the columns of an orthogonal array
will not change the  corresponding Latin square.

On the other hand, if we can construct a matrix $V$ = $\left[\begin{array}{ccccc}
a_{1} & a_{2} & a_{3} & \cdots\cdots\  & a_{n^{2}}\\
b_{1} & b_{2} & b_{3} & \cdots\cdots\  & b_{n^{2}}\\
c_{1} & c_{2} & c_{3} & \cdots\cdots\  & c_{n^{2}}
\end{array}\right]$ of size $3\times n^{2}$  by  satisfying the following conditions:
\label{Cndt:Orth-Arr} 
\begin{itemize}
\item (O1) Every row of the array is comprised  the elements of  1, 2, $\cdots$,
$n$; 
\item (O2) Every element $k$ (1$\leqslant k\leqslant n$) appears  exactly
 $n$ times in  each  row; 
\item (O3) The columns of this array are orthogonal pairwise, that is, any
pair of columns  shares  at most one element in the same position; 
\end{itemize}
then we can construct matrix $Y_{2}$ of order $n$ from this orthogonal
array by putting number $c_{t}$ in the position  of  $\left(a_{t},\:b_{t}\right)$
of an empty matrix of order $n$ ($t$=1, 2, 3, $\cdots$, $n^{2}$).
It is clear that $Y_{2}$ is a Latin square. Hence,  array $V$ is
the orthogonal array of Latin square $Y_{2}$.

So,  there is a  one-to-one  correspondence between the Latin squares
of order $n$ and the arrays of size $3\times n^{2}$ satisfying the
 earlier conditions  ( if the ordering of the columns of the orthogonal
array is ignored ). Hence,  there is no problem to call a $3\times n^{2}$
matrix  as  an orthogonal array if it satisfies the 3 conditions mentioned
here.

This means that an orthogonal array of size $3\times n^{2}$ corresponding
to Latin square $Y$ will still be an orthogonal array after permuting
its rows. But  the new orthogonal array will correspond to another
Latin square (called  as  a \emph{conjugate} or \emph{parastrophe},
which will be  stated  later) closely related to $Y$.

\label{lem:Iso_on_OA} The operation of permuting the rows of Latin
square $Y$ according to permutation $\alpha$ will correspond to
the action of $\alpha$ on the members in the 1st row of its orthogonal
array. Permuting the columns of Latin square $Y$ according to permutation
$\beta$ will correspond to replacing the number $i$ in the  second
 row of its orthogonal array  with  $\beta(i)$ ($i$ = 1, 2, $\cdots$,
$n$). The transformation $\mathscr{L}_{\gamma}$ acting on Latin
square $Y$ corresponds to the action of substituting the number $i$
 with  $\gamma(i)$ in the  third  row of the orthogonal array ($i$
= 1, 2, $\cdots$, $n$).

Here,  we define the representative of an isotopy class or a main
class as the minimal one in  the  lexicographic order. Obviously,
there are some other standards for an isotopy or a main class representative.
For some purpose, some other definition of  the  canonical form of
an equivalence class will be more efficient in  practice. These definitions
will not be discussed here, but the related algorithms derived from
the conclusions in this paper, may rely on them.

\section{parastrophes of a Latin square\label{subsec:Adjugates-LS}}

It  is not  difficult to understand that the permutation of  the  rows
of an orthogonal array will not interfere the orthogonality of  its
 columns. The array after row permutation will also remain as  an
orthogonal array of a certain Latin square. The new Latin square  will
be  firmly related to the original one. As the number of  reordering
 of a sequence of 3 different entries is 6, there are  exactly  6
transformations for permuting the rows of an orthogonal array.

Every column of an orthogonal array consists of 3 elements, the first
 one  is the row index, the second  one  is the column index, and
the third  one  is the entry of the Latin square in the position determined
by the  previous two elements. Let $[\mathrm{r},\mathrm{c},\mathrm{e}]$
or (1)  denotes  the transformation to keep the original order of
the rows of an orthogonal array; $[\mathrm{c},\mathrm{r},\mathrm{e}]$
or $(1\ 2)$  denotes  the transformation to interchange rows 1 and
2 of an orthogonal array;  and  $[\mathrm{r},\mathrm{e},\mathrm{c}]$
or $(2\ 3)$  denotes  the transformation to interchange rows 2 and
3 of an orthogonal array, etc.  A  similar notation may be found in
reference \cite{Brendan2007SmallLS}  also. There are 5 transformations
that will indeed change the order of the rows of an orthogonal array.
For convenience, denote the $i$'th row of Latin square $Y$ by $Y_{i}$
($i$ = 1, 2, $\cdots$, $n$). According to the convention  defined
on page \pageref{CVT:Permutation-2},  sometimes  $Y_{i}$ = $\left[\begin{array}{cccc}
y_{i1}, & y_{i2}, & \cdots, & y_{in}\end{array}\right]$ will represent the transformation that sends $j$ to $y_{ij}$ ($j$
= 1, 2, $\cdots$, $n$), and permutation $\alpha\in\mathrm{S}_{n}$
will denote the sequence {[}$\alpha(1),$ $\alpha(2)$, $\cdots$,
$\alpha(n)${]},  which  will rely  on the  corresponding  contexts.

 (1) $[\mathrm{r},\mathrm{e},\mathrm{c}]$ or $(2\ 3)$  \\
 If the  second and third rows  of an orthogonal array $V$ are interchanged,
it will result  to:  
\begin{equation}
V^{(\mathrm{I})}=\left[\begin{array}{llll|llll|lll|llll}
1 & 1 & \cdots & 1 & \,2 & 2 & \cdots & 2 & \,\bullet & \bullet & \bullet\, & \,n & n & \cdots & n\\
y_{11} & y_{12} & \cdots & y_{1n} & \,y_{21} & y_{22} & \cdots & y_{2n} & \,\bullet & \bullet & \bullet & \,y_{n1} & y_{n2} & \cdots & y_{nn}\\
1 & 2 & \cdots & n & \,1 & 2 & \cdots & n & \,\bullet & \bullet & \bullet & \,1 & 2 & \cdots & n
\end{array}\right].\label{eq:OrthogonalArray-2}
\end{equation}
Sort the columns of $V^{(\mathrm{I})}$ in  the  lexicographical order,
such that the  sub-matrix consisted  of the  first and second rows
 will be :  
\begin{equation}
V_{0}=\left[\begin{array}{cccc|cccc|ccc|cccc}
1 & 1 & \cdots & 1\, & \,2 & 2 & \cdots & 2\, & \,\bullet & \bullet & \bullet\, & \,n & n & \cdots & n\\
1 & 2 & \cdots & n\, & \,1 & 2 & \cdots & n\, & \,\bullet & \bullet & \bullet\, & \,1 & 2 & \cdots & n
\end{array}\right].\label{eq:OrthogonalArray_2Row}
\end{equation}
It results  to:  
\begin{equation}
V^{(\mathrm{IA})}=\left[\begin{array}{llll|llll|lll|llll}
1 & 1 & \cdots & 1 & \,2 & 2 & \cdots & 2 & \,\bullet & \bullet & \bullet\, & \,n & n & \cdots & n\\
1 & 2 & \cdots & n & \,1 & 2 & \cdots & n & \,\bullet & \bullet & \bullet & \,1 & 2 & \cdots & n\\
y'_{11} & y'_{12} & \cdots & y'_{1n} & \,y'_{21} & y'_{22} & \cdots & y'_{2n} & \,\bullet & \bullet & \bullet & \,y'_{n1} & y'_{n2} & \cdots & y'_{nn}
\end{array}\right].\label{eq:OrthogonalArray-2b}
\end{equation}
 Here,  $V^{(\mathrm{IA})}$ is essentially the same as $V^{(\mathrm{I})}$
as they  both  correspond to the same Latin square. There are  ($n-1$)
 vertical lines in  array  $V$, which divide it into $n$ segments.
Every segment of $V$ represents a row of $Y$ as the members in a
segment have the same row index. It is clear that the columns of $V^{(\mathrm{I})}$
in any segment will  remain  in that segment after  the  sorting  also.
That is, the entries  of a  row of the original Latin square $Y$
will still be in the same row in the new Latin square (denoted by
$Y^{(\mathrm{I})}$) \label{Sym:Y^(I)}  corresponding  to the orthogonal
array $V^{(\mathrm{I})}$  as  the row index of every entry is not
changed in the process of interchanging the  second and third rows
 of the orthogonal array. In every segment, the  second and third
rows  will  create  a permutation in  the  two-row form. It is widely
known that  the exchanging of  the two rows of permutation $\alpha$=$\left(\begin{array}{cccc}
1 & 2 & \cdots & n\\
a_{1} & a_{2} & \cdots & a_{n}
\end{array}\right)$ will result  to  its inverse $\alpha^{-1}$ = $\left(\begin{array}{cccc}
a_{1} & a_{2} & \cdots & a_{n}\\
1 & 2 & \cdots & n
\end{array}\right)$ = $\left(\begin{array}{cccc}
1 & 2 & \cdots & n\\
a'_{1} & a'_{2} & \cdots & a'_{n}
\end{array}\right)$. So,  the sequence $\left[\begin{array}{cccc}
y'_{i1}, & y'_{i2}, & \cdots, & y'_{in}\end{array}\right]$ is the inverse of $\left[\begin{array}{cccc}
y_{i1}, & y_{i2}, & \cdots, & y_{in}\end{array}\right]$ as both are  the reordering  of $\left[\begin{array}{cccc}
1, & 2, & \cdots, & n\end{array}\right]$ ($i$ = 1, 2, $\cdots$, $n$). Hence,  the interchanging of the
second and third rows  of the orthogonal array of a Latin square corresponds
to  substitution of  every row of the Latin square by its inverse,
i.e., $Y^{(\mathrm{I})}$ = $\left(\begin{array}{c}
Y_{1}^{-1}\\
\vdots\\
Y_{n}^{-1}
\end{array}\right)$. This operation is mentioned implicitly in reference \cite{Lam1990NumLS8}.

 (2) $[\mathrm{c},\mathrm{r},\mathrm{e}]$ or $(1\ 2)$ \\
 With regard to the orthogonal array $V$ in \eqref{eq:OrthogonalArray}
of Latin square $Y=\left(y_{ij}\right)_{n\times n}$,  the interchanging
of the first and second rows of $V$ will result to:  
\begin{equation}
V^{(\mathrm{II})}=\left[\begin{array}{llll|llll|lll|llll}
1 & 2 & \cdots & n & \,1 & 2 & \cdots & n & \,\bullet & \bullet & \bullet & \,1 & 2 & \cdots & n\\
1 & 1 & \cdots & 1 & \,2 & 2 & \cdots & 2 & \,\bullet & \bullet & \bullet\, & \,n & n & \cdots & n\\
y_{11} & y_{12} & \cdots & y_{1n} & \,y_{21} & y_{22} & \cdots & y_{2n} & \,\bullet & \bullet & \bullet & \,y_{n1} & y_{n2} & \cdots & y_{nn}
\end{array}\right].\label{eq:OrthogonalArray-1}
\end{equation}
It means that every entry $(i,j,y_{ij})^{\mathrm{T}}$ will become
$(j,i,y_{ij})^{\mathrm{T}}$. That is, the entry in the $(j,i)$ -th
 position of the new Latin  square  (denoted by $Y^{(\mathrm{II})}$)
 corresponding  to $V^{(\mathrm{II})}$ is the entry $y_{ij}$ in
the $(i,j)$ -th  position of the original Latin square $Y$. Hence,
 the new Latin square $Y^{(\mathrm{II})}$ is the transpose of the
original one.

In order to understand the transformation $[\mathrm{c},\mathrm{e},\mathrm{r}]$,
we should first introduce the operation $[\mathrm{e},\mathrm{c},\mathrm{r}]$.

 (3) $[\mathrm{e},\mathrm{r},\mathrm{c}]$ or $(1\ 3\ 2)$  \\
  Firstly, interchange rows 1 and 2 of $V$,  and  then interchange
rows 2 and 3. This operation corresponds to replacing the $i$'th
row of $Y$ with the inverse of the $i$'th row of its transpose $Y^{\mathrm{T}}$,
or substituting the $i$'th row of $Y$ by the inverse of the $i$'th
column of $Y$.  Express  the result  as  $Y^{(\mathrm{III})}$ =
$\left(\begin{array}{c}
\left(Z_{1}^{-1}\right)^{\mathrm{T}}\\
\vdots\\
\left(Z_{n}^{-1}\right)^{\mathrm{T}}
\end{array}\right)$. Here,  the  superscript  ``$\mathrm{T}$'' means the transpose
 as  $Z_{i}$ is a column of $Y$ ($i$ = 1, 2, $\cdots$, $n$).

 (4) $[\mathrm{e},\mathrm{c},\mathrm{r}]$ or $(1\ 3)$ \\
 If the  first and third rows  of $V$ are interchanged, it will result
 to:  
\begin{equation}
V^{(\mathrm{IV})}=\left[\begin{array}{llll|llll|lll|llll}
y_{11} & y_{12} & \cdots & y_{1n} & \,y_{21} & y_{22} & \cdots & y_{2n} & \,\bullet & \bullet & \bullet & \,y_{n1} & y_{n2} & \cdots & y_{nn}\\
1 & 2 & \cdots & n & \,1 & 2 & \cdots & n & \,\bullet & \bullet & \bullet & \,1 & 2 & \cdots & n\\
1 & 1 & \cdots & 1 & \,2 & 2 & \cdots & 2 & \,\bullet & \bullet & \bullet\, & \,n & n & \cdots & n
\end{array}\right].\label{eq:OrthogonalArray-3}
\end{equation}
It is difficult to find the  relationship between  the original Latin
square $Y$ and the new Latin square $V^{(\mathrm{IV})}$  corresponding
 to array $V^{(\mathrm{IV})}$ by sorting the columns of $V^{(\mathrm{IV})}$
in  the  lexicographic order. But  it will work in another way also.
 The interchanging of columns $(l-1)n+k$ and $(k-1)n+l$ of $V$,
($l$ = 1, 2, $\cdots$, $n-1$; $k$ = $l+1$, $\cdots$, $n$),
means the sorting of  the columns of $V$  so  that the  sub-matrix
consisted of the first and second rows becomes:  
\begin{equation}
V_{0}^{(\mathrm{A})}=\left[\begin{array}{cccc|cccc|ccc|cccc}
1 & 2 & \cdots & n\, & \,1 & 2 & \cdots & n\, & \,\bullet & \bullet & \bullet\, & \,1 & 2 & \cdots & n\\
1 & 1 & \cdots & 1\, & \,2 & 2 & \cdots & 2\, & \,\bullet & \bullet & \bullet\, & \,n & n & \cdots & n
\end{array}\right],\label{eq:OrthogonalArray_2Row-B}
\end{equation}
 which will result to:  
\begin{equation}
V^{(\mathrm{A})}=\left[\begin{array}{llll|llll|lll|llll}
1 & 2 & \cdots & n & \,1 & 2 & \cdots & n & \,\bullet & \bullet & \bullet & \,1 & 2 & \cdots & n\\
1 & 1 & \cdots & 1 & \,2 & 2 & \cdots & 2 & \,\bullet & \bullet & \bullet\, & \,n & n & \cdots & n\\
z_{11} & z_{12} & \cdots & z_{1n} & \,z_{21} & z_{22} & \cdots & z_{2n} & \,\bullet & \bullet & \bullet & \,z_{n1} & z_{n2} & \cdots & z_{nn}
\end{array}\right].\label{eq:OrthogonalArray-b}
\end{equation}
 Here also, $V^{(\mathrm{A})}$ corresponds to Latin square $Y$.
The vertical lines in the array $V^{(\mathrm{A})}$ divide it into
$n$ segments. Every segment corresponds to a column of $Y$  as  all
the entries in a segment of $V^{(\mathrm{A})}$ share the same column
index. That is, {[} $z_{i1}$, $z_{i2}$, $\cdots$, $z_{in}$ {]}
is the $i$'th column of $Y$.

 The interchanging of the first and third rows  of $V^{(\mathrm{A})}$
will result  to:  
\begin{equation}
V^{(\mathrm{IVA})}=\left[\begin{array}{llll|llll|lll|llll}
z_{11} & z_{12} & \cdots & z_{1n} & \,z_{21} & z_{22} & \cdots & z_{2n} & \,\bullet & \bullet & \bullet & \,z_{n1} & z_{n2} & \cdots & z_{nn}\\
1 & 1 & \cdots & 1 & \,2 & 2 & \cdots & 2 & \,\bullet & \bullet & \bullet\, & \,n & n & \cdots & n\\
1 & 2 & \cdots & n & \,1 & 2 & \cdots & n & \,\bullet & \bullet & \bullet & \,1 & 2 & \cdots & n
\end{array}\right].\label{eq:OrthogonalArray-3-A}
\end{equation}
Obviously, $V^{(\mathrm{IVA})}$ and $V^{(\mathrm{IV})}$ correspond
to the same Latin square. Sorting the columns in each segment of $V^{(\mathrm{IVA})}$,
 such that the  first and second rows  are the same as  those of  $V_{0}^{(\mathrm{A})}$,
 the following is obtained:  
\begin{equation}
V^{(\mathrm{IVB})}=\left[\begin{array}{llll|llll|lll|llll}
1 & 2 & \cdots & n & \,1 & 2 & \cdots & n & \,\bullet & \bullet & \bullet & \,1 & 2 & \cdots & n\\
1 & 1 & \cdots & 1 & \,2 & 2 & \cdots & 2 & \,\bullet & \bullet & \bullet\, & \,n & n & \cdots & n\\
z'_{11} & z'_{12} & \cdots & z'_{1n} & \,z'_{21} & z'_{22} & \cdots & z'_{2n} & \,\bullet & \bullet & \bullet & \,z'_{n1} & z'_{n2} & \cdots & z'_{nn}
\end{array}\right].\label{eq:OrthogonalArray-3-B}
\end{equation}

So,  {[} $z'_{i1}$, $z'_{i2}$, $\cdots$, $z'_{in}$ {]} is the
inverse of {[} $z_{i1}$, $z_{i2}$, $\cdots$, $z_{in}$ {]}. Let
the Latin square  corresponding  to $V^{(\mathrm{IVB})}$ (or $V^{(\mathrm{IVA})}$,
$V^{(\mathrm{IV})}$) be $Y^{(\mathrm{IV})}$. Therefore, the Latin
square $Y^{(\mathrm{IV})}$,  related to $V^{(\mathrm{IV})}$  and
 generated by  interchanging  the  first and third rows  of the orthogonal
array $V$, consists of the columns,  which are the inverse of the
columns of $Y$.  $Y^{(\mathrm{IV})}$ = $\left(\begin{array}{ccc}
Z_{1}^{-1}\  & \cdots\  & Z_{n}^{-1}\end{array}\right)$. 

 (5) $[\mathrm{c},\mathrm{e},\mathrm{r}]$ or $(1\ 2\ 3)$  \\
  First  interchange rows 1 and 3 of $V$,  and  then interchange
rows 2 and 3, or equivalently, first interchange rows 1 and 2,  and
 then interchange rows 1 and 3. This operation corresponds to substituting
the $i$'th column of $Y$ by the inverse of the $i$'th column of
its transpose $Y^{\mathrm{T}}$, or replacing the $i$'th column of
$Y$ with the inverse of the $i$'th row of $Y$.  Express  the result
 as  $Y^{(\mathrm{V})}$. \label{Sym:Y^(V)} $Y^{(V)}$ = $\left(\begin{array}{ccc}
\left(Y_{1}^{-1}\right)^{\mathrm{T}}\  & \cdots\  & \left(Y_{n}^{-1}\right)^{\mathrm{T}}\end{array}\right)$.

$\ $
\begin{lem}
\label{lem:Rel-Ad-T} Let $\rho\in\mathrm{S}_{3}$,  denote by $\mathcal{F}_{\rho}$
the transformation mentioned above, i.e., $\mathcal{F}_{(1)}(Y)$
= $Y$, $\mathcal{F}_{(1\,2)}(Y)$ = $Y^{(\mathrm{II})}$ = $Y^{\mathrm{T}}$,
$\mathcal{F}_{(2\,3)}(Y)$ = $Y^{(\mathrm{I})}$, $\mathcal{F}_{(1\,3)}(Y)$
= $Y^{(\mathrm{IV})}$, etc. By definition, it is clear that :  
\begin{equation}
\mathcal{F}_{\rho_{1}}\circ\mathcal{F}_{\rho_{2}}=\mathcal{F}_{\rho_{1}\rho_{2}}\quad(\forall\rho_{1},\rho_{2}\in\mathrm{S}_{3}).
\end{equation}
\end{lem}

The 6 Latin squares , corresponding  to the orthogonal arrays obtained
by permuting the rows of the orthogonal array $V$ of the Latin square
$Y$,  are called the \emph{conjugate}s\index{conjugates (of a Latin square)}
or \emph{adjugates\index{adjugate}} or \emph{parastrophes} \index{parastrophe}
of $Y$. Of course, $Y$ is a conjugate of itself. The set of the
Latin squares that are isotopic to any conjugate of Latin square $Y$
is called the \emph{main class}\index{main class} or \emph{specy
\index{specy}} or \emph{paratopy class \index{paratopy class}} of
$Y$. If Latin square $Z$ belongs to the main class of $Y$, then
$Y$ and $Z$ are called  the  \emph{paratopic} \index{paratopic}
or \emph{main class equivalent\index{main class equivalent}}. Sometimes,
the set of the Latin squares that are isotopic to $Y$ or $Y^{\mathrm{T}}$
is called the \emph{type} \index{type} of $Y$  (refer to \cite{Bailey2003Equivalence}
or \cite{Brendan2007SmallLS}). 

In this paper, a special equivalence class, namely the \emph{inverse
type\index{inverse type}},  is  defined for  the convenience in the
process of generating the representatives of all the main classes
of Latin squares of a certain order. The \emph{inverse type} of Latin
square $Y$ is the set of the Latin squares that are isotopic to $Y$
or $Y^{(\mathrm{I})}$, where $Y^{(\mathrm{I})}$ is the Latin square
 corresponding  to the orthogonal array $V^{(\mathrm{I})}$ obtained
by interchanging the  second and third rows of  orthogonal array $V$
of $Y$ as described before  (this idea is adopted  from the notion
 of  ``row inverse”  defined  in reference \cite{Lam1990NumLS8}).
 The \emph{inverse type} mentioned here may be called  as the  ``row
inverse type'' for the sake of accuracy  as  $Y^{(\mathrm{IV})}$
is another type of inverse (column inverse) of $Y$. The reason for
choosing  the  row inverse type is that the Latin squares are generated
by rows  in the following papers by the authors   (in  some papers,
 the  Latin squares are generated by columns,  so that  the column
inverse type will be useful). 

\section{Relations on Paratopisms\label{subsec:RLT-Prtp}}

The transformation that  sends  a Latin square to  its another paratopic
 is called a \emph{paratopism}\index{paratopism} or a \emph{paratopic
transformation.}\index{paratopic transformation} Let $\mathcal{P}_{n}$
\label{Sym:mathcal_P_n} \nomenclature[P_n]{ $\mathcal{P}_{n}$ or  $\mathcal{P}$}{ the set of all the paratopic transformations of Latin squares of order $n$. \pageref{Sym:mathcal_P_n}}
be the set of all the paratopic transformations of Latin squares of
order $n$, sometimes denoted by $\mathcal{P}$  in  short if it results
no ambiguity. It is obvious that $\mathcal{P}$ together with the
composition operation ``$\circ$'' will form a group, called the
\emph{paratopic transformation group} \index{paratopic transformation group (of Latin squares)}
or \emph{paratopism} \emph{group}\index{paratopism group@paratopism\emph{ }group}.
Obviously, all the isotopy transformations are paratopisms.

It will not be difficult to find out by hand  calculation  the following
relations of $\mathcal{F}_{\rho}$  with  $\mathscr{C}_{\alpha}$,
$\mathscr{R}_{\beta}$  and  $\mathscr{L}_{\gamma}$ when  their orders
are exchanged  if we are familiar with how the transformations $\mathscr{C}_{\alpha}$,
$\mathscr{R}_{\beta}$,  and  $\mathscr{L}_{\gamma}$ change the rows
and columns in detail. 
\begin{thm}
\label{thm:Cjgt-ISTP-Comsition-order-1}For $\forall\mathscr{T}\in\mathscr{I}_{n}$,
$\forall\alpha,\beta,\gamma\in\mathscr{\mathrm{S}}_{n}$,  the following
equalities hold: \\
 \hphantom{$\ MM$ } $\mathscr{T}\circ\mathcal{F}_{(1)}$ = $\mathcal{F}_{(1)}\circ\mathscr{T}$,
\\
 \hphantom{$\ MM$ } $\mathscr{R}_{\alpha}\circ\mathcal{F}_{(1\,2)}$=$\mathcal{F}_{(1\,2)}\circ\mathscr{C}_{\alpha}$,
$\quad\ \mathscr{C}_{\beta}\circ\mathcal{F}_{(1\,2)}$=$\mathcal{F}_{(1\,2)}\circ\mathscr{R}_{\beta}$,
$\mathscr{\quad\ L}_{\gamma}\circ\mathcal{F}_{(1\,2)}$=$\mathcal{F}_{(1\,2)}\circ\mathscr{L}_{\gamma}$,\\
 \hphantom{$\ MM$ } $\mathscr{R}_{\alpha}\circ\mathcal{F}_{(2\,3)}$=$\mathcal{F}_{(2\,3)}\circ\mathscr{R}_{\alpha}$,$\quad\ \mathscr{C}_{\beta}\circ\mathcal{F}_{(2\,3)}$=$\mathcal{F}_{(2\,3)}\circ\mathscr{L}_{\beta}$,
$\mathscr{\quad\ L}_{\gamma}\circ\mathcal{F}_{(2\,3)}$=$\mathcal{F}_{(2\,3)}\circ\mathscr{C}_{\gamma}$,\\
 \hphantom{$\ MM$ } $\mathscr{R}_{\alpha}\circ\mathcal{F}_{(1\,3)}$=$\mathcal{F}_{(1\,3)}\circ\mathscr{L}_{\alpha}$,
$\quad\ \mathscr{C}_{\beta}\circ\mathcal{F}_{(1\,3)}$=$\mathcal{F}_{(1\,3)}\circ\mathscr{C}_{\beta}$,
$\mathscr{\quad\ L}_{\gamma}\circ\mathcal{F}_{(1\,3)}$=$\mathcal{F}_{(1\,3)}\circ\mathscr{R}_{\gamma}$. 
\end{thm}

The equalities in the  first and second lines  are obvious. Here,
 we explain some equalities in the  third and fourth lines. The readers
may get the proof of other equalities in the same way without  any
difficulty.

 As per the convention described  on page \pageref{CVT:Permutation-3}
: 

$Y=\left(\begin{array}{c}
Y_{1}\\
\vdots\\
Y_{n}
\end{array}\right)$ $\overset{\mathcal{F}_{(2\,3)}}{\longrightarrow}$ $\left(\begin{array}{c}
Y_{1}^{-1}\\
\vdots\\
Y_{n}^{-1}
\end{array}\right)$ $\overset{\mathscr{R}_{\alpha}}{\longrightarrow}$ $\left(\begin{array}{c}
Y_{\alpha^{-1}(1)}^{-1}\\
\vdots\\
Y_{\alpha^{-1}(n)}^{-1}
\end{array}\right)$ $\overset{\mathcal{F}_{(2\,3)}}{\longleftarrow}$ $\left(\begin{array}{c}
Y_{\alpha^{-1}(1)}\\
\vdots\\
Y_{\alpha^{-1}(n)}
\end{array}\right)$ $\overset{\mathscr{R}_{\alpha}}{\longleftarrow}$ $\left(\begin{array}{c}
Y_{1}\\
\vdots\\
Y_{n}
\end{array}\right)$, hence $\mathscr{R}_{\alpha}\circ\mathcal{F}_{(2\,3)}$=$\mathcal{F}_{(2\,3)}\circ\mathscr{R}_{\alpha}$
holds.

$Y=\left(\begin{array}{c}
Y_{1}\\
\vdots\\
Y_{n}
\end{array}\right)$ $\overset{\mathcal{F}_{(2\,3)}}{\longrightarrow}$ $\left(\begin{array}{c}
Y_{1}^{-1}\\
\vdots\\
Y_{n}^{-1}
\end{array}\right)$ $\overset{\mathscr{C}_{\beta}}{\longrightarrow}$ $\left(\begin{array}{c}
Y_{1}^{-1}\beta^{-1}\\
\vdots\\
Y_{n}^{-1}\beta^{-1}
\end{array}\right)$ $\overset{\mathcal{F}_{(2\,3)}}{\longleftarrow}$ $\left(\begin{array}{c}
\beta Y_{1}\\
\vdots\\
\beta Y_{n}
\end{array}\right)$ $\overset{\mathscr{L}_{\beta}}{\longleftarrow}$ $\left(\begin{array}{c}
Y_{1}\\
\vdots\\
Y_{n}
\end{array}\right)$, so $\ $ \\
 $\mathscr{C}_{\beta}\circ\mathcal{F}_{(2\,3)}$ = $\mathcal{F}_{(2\,3)}\circ\mathscr{L}_{\beta}$.

$Y=\left(\begin{array}{c}
Y_{1}\\
\vdots\\
Y_{n}
\end{array}\right)$ $\overset{\mathcal{F}_{(2\,3)}}{\longrightarrow}$ $\left(\begin{array}{c}
Y_{1}^{-1}\\
\vdots\\
Y_{n}^{-1}
\end{array}\right)$ $\overset{\mathscr{L}_{\gamma}}{\longrightarrow}$ $\left(\begin{array}{c}
\gamma Y_{1}^{-1}\\
\vdots\\
\gamma Y_{n}^{-1}
\end{array}\right)$ $\overset{\mathcal{F}_{(2\,3)}}{\longleftarrow}$ $\left(\begin{array}{c}
Y_{1}\gamma^{-1}\\
\vdots\\
Y_{n}\gamma^{-1}
\end{array}\right)$ $\overset{\mathscr{C}_{\gamma}}{\longleftarrow}$ $\left(\begin{array}{c}
Y_{1}\\
\vdots\\
Y_{n}
\end{array}\right)$, therefore $\mathscr{L}_{\gamma}\circ\mathcal{F}_{(2\,3)}$=$\mathcal{F}_{(2\,3)}\circ\mathscr{C}_{\gamma}$.

$Y$=$\left(\begin{array}{ccc}
Z_{1}, & \cdots, & Z_{n}\end{array}\right)$ $\overset{\mathcal{F}_{(1\,3)}}{\longrightarrow}$ $\left(\begin{array}{ccc}
Z_{1}^{-1}, & \cdots, & Z_{n}^{-1}\end{array}\right)$ $\overset{\mathscr{R}_{\alpha}}{\longrightarrow}$ $\left(\begin{array}{ccc}
Z_{1}^{-1}\alpha^{-1}, & \cdots, & Z_{n}^{-1}\alpha^{-1}\end{array}\right)$,\\
 $Y$=$\left(\begin{array}{ccc}
Z_{1}, & \cdots, & Z_{n}\end{array}\right)$ $\overset{\mathscr{L}_{\alpha}}{\longrightarrow}$ $\left(\begin{array}{ccc}
\alpha Z_{1}, & \cdots, & \alpha Z_{n}\end{array}\right)$ $\overset{\mathcal{F}_{(1\,3)}}{\longrightarrow}$ $\left(\begin{array}{ccc}
Z_{1}^{-1}\alpha^{-1}, & \cdots, & Z_{n}^{-1}\alpha^{-1}\end{array}\right)$,\\
 then $\mathscr{R}_{\alpha}\circ\mathcal{F}_{(1\,3)}$=$\mathcal{F}_{(1\,3)}\circ\mathscr{L}_{\alpha}$.

With the  above equalities, it will be easy to obtain the properties
when $\mathcal{F}_{(1\,3\,2)}$ or $\mathcal{F}_{(1\,2\,3)}$  combined
 with $\mathscr{R}_{\alpha}$, $\mathscr{C}_{\beta}$, or $\mathscr{L}_{\gamma}$.
 Then,  we will have  the following:  
\begin{thm}
\label{thm:Cjgt-ISTP-Comsition-order-2}For $\forall\alpha,\beta,\gamma\in\mathscr{\mathrm{S}}_{n}$,\textup{
}$\forall\mathscr{T}\in\mathscr{I}_{n}$,\\
 \hphantom{$\ MM$ }$\mathscr{T}\circ\mathcal{F}_{(1)}=\mathcal{F}_{(1)}\circ\mathscr{T}$,\\
 \hphantom{$\ MM$ }\textup{$\left(\mathscr{R}_{\alpha}\circ\mathscr{C}_{\beta}\circ\mathscr{L}_{\gamma}\right)\circ\mathcal{F}_{(1\,2)}=\mathcal{F}_{(1\,2)}\circ\left(\mathscr{R}_{\beta}\circ\mathscr{C}_{\alpha}\circ\mathscr{L}_{\gamma}\right)$,
}\\
 \hphantom{$\ MM$ }\textup{$\left(\mathscr{R}_{\alpha}\circ\mathscr{C}_{\beta}\circ\mathscr{L}_{\gamma}\right)\circ\mathcal{F}_{(1\,3)}=\mathcal{F}_{(1\,3)}\circ\left(\mathscr{R}_{\gamma}\circ\mathscr{C}_{\beta}\circ\mathscr{L}_{\alpha}\right)$,}\\
 \hphantom{$\ MM$ }\textup{$\left(\mathscr{R}_{\alpha}\circ\mathscr{C}_{\beta}\circ\mathscr{L}_{\gamma}\right)\circ\mathcal{F}_{(2\,3)}=\mathcal{F}_{(2\,3)}\circ\left(\mathscr{R}_{\alpha}\circ\mathscr{C}_{\gamma}\circ\mathscr{L}_{\beta}\right)$.} 
\end{thm}

Since $\mathcal{F}_{(1\,2\,3)}$ =$\mathcal{F}_{(1\,3)}$$\circ\mathcal{F}_{(1\,2)}$,
$\mathcal{F}_{(1\,3\,2)}$ =$\mathcal{F}_{(1\,2)}$$\circ\mathcal{F}_{(1\,3)}$,
so 
\begin{thm}
\label{thm:Cjgt-ISTP-Comsition-order-3}For $\forall\alpha,\beta,\gamma\in\mathscr{\mathrm{S}}_{n}$,\\
 \hphantom{$\ MM$ }\textup{$\left(\mathscr{R}_{\alpha}\circ\mathscr{C}_{\beta}\circ\mathscr{L}_{\gamma}\right)\circ\mathcal{F}_{(1\,2\,3)}=\mathcal{F}_{(1\,2\,3)}\circ\left(\mathscr{R}_{\beta}\circ\mathscr{C}_{\gamma}\circ\mathscr{L}_{\alpha}\right)$,
}\\
 \hphantom{$\ MM$ }\textup{$\left(\mathscr{R}_{\alpha}\circ\mathscr{C}_{\beta}\circ\mathscr{L}_{\gamma}\right)\circ\mathcal{F}_{(1\,3\,2)}=\mathcal{F}_{(1\,3\,2)}\circ\left(\mathscr{R}_{\gamma}\circ\mathscr{C}_{\alpha}\circ\mathscr{L}_{\beta}\right)$.
} 
\end{thm}

Here,  we denote an isotopism by $\mathscr{R}_{\alpha}\circ\mathscr{C}_{\beta}\circ\mathscr{L}_{\gamma}$,
not by  $\mathscr{C}_{\beta}\circ\mathscr{R}_{\alpha}\circ\mathscr{L}_{\gamma}$
(although $\mathscr{R}_{\alpha}$, $\mathscr{C}_{\beta}$, $\mathscr{L}_{\gamma}$
commute pairwise). The reason is that,  when exchanging the position
of $\left(\mathscr{R}_{\alpha}\circ\mathscr{C}_{\beta}\circ\mathscr{L}_{\gamma}\right)$
and $\mathcal{F}_{\rho}$,  the subscripts of the three transformations
are permuted  according to the permutation $\rho\in\mathscr{\mathrm{S}}_{3}$
 as  shown above  (just substitute 1, 2, 3 by $\alpha,\beta,\gamma$,
respectively. For instance, $(1\,2\,3)$ will become $\left(\alpha\,\beta\,\gamma\right)$.
So, when moving $\mathcal{F}_{(1\,2\,3)}$ from the right side of
$\left(\mathscr{R}_{\alpha}\circ\mathscr{C}_{\beta}\circ\mathscr{L}_{\gamma}\right)$
to the left side, $\alpha,\beta,\gamma$ will become $\beta,\gamma,\alpha$,
respectively).  So,  we can denote these formulas as below.
\begin{thm}
\emph{(main result 1)} \label{thm:Cjgt-ISTP-Comsition-order-General}
$\forall\,\alpha_{1},\alpha_{2},\alpha_{3}\in\mathscr{\mathrm{S}}_{n}$,
\textup{$\forall\,\beta_{1},\beta_{2},\beta_{3}\in\mathscr{\mathrm{S}}_{n}$,
}$\forall\rho\in\mathscr{\mathrm{S}}_{3}$, 
\begin{align}
\left(\mathscr{R}_{\alpha_{1}}\circ\mathscr{C}_{\alpha_{2}}\circ\mathscr{L}_{\alpha_{3}}\right)\circ\mathcal{F}_{\rho} & =\mathcal{F}_{\rho}\circ\left(\mathscr{R}_{\alpha_{\rho(1)}}\circ\mathscr{C}_{\alpha_{\rho(2)}}\circ\mathscr{L}_{\alpha_{\rho(3)}}\right),\label{eq:Cjgt-ISTP-Comsition-order-1}\\
\mathcal{F}_{\rho}\circ\left(\mathscr{R}_{\beta_{1}}\circ\mathscr{C}_{\beta_{2}}\circ\mathscr{L}_{\beta_{3}}\right) & =\left(\mathscr{R}_{\beta_{\rho^{-1}(1)}}\circ\mathscr{C}_{\beta_{\rho^{-1}(2)}}\circ\mathscr{L}_{\beta_{\rho^{-1}(3)}}\right)\circ\mathcal{F}_{\rho}.\label{eq:Cjgt-ISTP-Comsition-order-2}
\end{align}
\end{thm}

With Lemmas  \ref{lem:Iso_on_OA} and \ref{lem:Rel-Ad-T}, it is not
difficult to explain the six theorems described above, although not
 very  intuitional.

It is clear that any paratopic transformation is the composition of
an isotopism $\left(\mathscr{R}_{\alpha}\circ\mathscr{C}_{\beta}\circ\mathscr{L}_{\gamma}\right)$
($\alpha,\beta,\gamma$ $\in\mathscr{\mathrm{S}}_{n}$) and certain
conjugate transformation $\mathcal{F}_{\rho}$ ($\rho\in\mathscr{\mathrm{S}}_{3}$).

Since $\mathscr{I}_{n}$=$\left\{ \left.\mathscr{R}_{\alpha}\circ\mathscr{C}_{\beta}\circ\mathscr{L}_{\gamma}\,\right|\,\alpha,\beta,\gamma\in\mathscr{\mathrm{S}}_{n}\,\right\} $
$\simeq$ $\mathscr{\mathrm{S}}_{n}^{3}$, it is convenient to denote
$\mathscr{R}_{\alpha}\circ\mathscr{C}_{\beta}\circ\mathscr{L}_{\gamma}\circ\mathcal{F}_{\rho}$
by $\mathscr{P}\left(\alpha,\beta,\gamma,\rho\right)$  in  short.
So,  we have $\left|\mathcal{P}_{n}\right|$ = $\left|\mathscr{\mathrm{S}}_{n}^{3}\right|\times\left|\mathscr{\mathrm{S}}_{3}\right|$
= $\left|\mathscr{I}_{n}\right|$$\cdot$$\left|\mathscr{\mathrm{S}}_{3}\right|$
= $6\left(n!\right)^{3}$. Let $\mathscr{P}$ : $\mathscr{\mathrm{S}}_{n}^{3}\times\mathscr{\mathrm{S}}_{3}$
$\rightarrow$ $\mathcal{P}_{n}$, $\left(\alpha,\beta,\gamma,\rho\right)$
$\longmapsto$ $\mathscr{R}_{\alpha}\circ\mathscr{C}_{\beta}\circ\mathscr{L}_{\gamma}\circ\mathcal{F}_{\rho}$.
It  is clear that $\mathscr{P}$ is a bijection. As a set, $\mathcal{P}_{n}$
is isomorphic to $\mathscr{I}_{n}$$\times$$\mathscr{\mathrm{S}}_{3}$
or $\mathscr{\mathrm{S}}_{n}^{3}\times\mathscr{\mathrm{S}}_{3}$.
But  as a group, $\mathcal{P}_{n}$ is not isomorphic to $\mathscr{I}_{n}$$\times$$\mathscr{\mathrm{S}}_{3}$
 (refer to \cite{Brendan2007SmallLS}),  because $\mathscr{P}$ is
not compatible with the multiplications in $\mathcal{P}_{n}$ as $\left(\mathscr{R}_{\alpha}\circ\mathscr{C}_{\beta}\circ\mathscr{L}_{\gamma}\right)$
and $\mathcal{F}_{\rho}$ do not commute unless $\rho=(1)$ or $\alpha=\beta=\gamma=(1)$.
\begin{thm}
\emph{(main result 2)} In general, $\forall\alpha_{1},\alpha_{2},\alpha_{3},$
$\beta_{1},\beta_{2},\beta_{3}$$\in\mathscr{\mathrm{S}}_{n}$, $\forall\rho,\zeta\in\mathscr{\mathrm{S}}_{3}$,

\hphantom{$\ MM$ } $\text{\textipa{\ }}$$\mathscr{P}\left(\alpha_{1},\alpha_{2},\alpha_{3},\rho\right)$
$\circ$ $\mathscr{P}\left(\beta_{1},\beta_{2},\beta_{3},\zeta\right)$
$=$ $\mathscr{P}\left(\alpha_{1}\beta_{\rho^{-1}(1)},\alpha_{2}\beta_{\rho^{-1}(2)},\alpha_{3}\beta_{\rho^{-1}(3)},\rho\zeta\right).$ 
\end{thm}

Usually, $\mathscr{P}\left(\alpha_{1}\beta_{\rho^{-1}(1)},\alpha_{2}\beta_{\rho^{-1}(2)},\alpha_{3}\beta_{\rho^{-1}(3)},\rho\zeta\right)$
differs from $\mathscr{P}\left(\alpha_{1}\beta_{1},\alpha_{2}\beta_{2},\alpha_{3}\beta_{3},\rho\zeta\right)$.

\section{Application}

 Using the  theorems mentioned above, we can avoid generating the
orthogonal array when producing the parastrophes of a Latin square,
which  improve the computational efficiency, especially  when writing
source codes.  Also, a lot of time can be saved in  generating all
the representatives of  the  main classes of Latin squares of a certain
order if we  use  these relations together with some properties of
cycle structures and properties of isotopic representatives, as we
can avoid generating a lot of Latin squares when testing  the  main
class representatives.

When generating the invariant group of a Latin square in paratopic
transformations, it will be more convenient to simplify the composition
of two or more paratopisms by  using  the relations described in Sec.\ \ref{lem:Rel-Ad-T}
(although this benefit  in improving the computational efficiency
is not very  conspicuous for a single Latin square, but it is remarkable
for a large number of Latin squares).

Some new algorithms for related problems, such as the algorithms for
testing or generating a representative of an equivalence class (a
main class or an isotopic class),  the algorithms for generating the
invariant group of a Latin square in some transformation (isotopisms
and main class transformations),  will be presented in near future.
 Those will be  different from  the ones  described in \cite{Petteri2006CACD}.

\section*{Acknowledgements}

\addcontentsline{toc}{section}{Acknowledgements}

The  authors  would like to express  their sincere  gratitude to Prof.\ \emph{LI
Shangzhi }from BUAA and Dr.\ ZHANG Zhe from Xidian University for
their valuable  suggestions for improving  this paper. The main part
of this paper is contained in the  Ph.D.  thesis \cite{liwenwei2014-LR-PP-RT-EATGC}
 of the first  author. The  reviewers, Professors Leonid A. Bokut,
Nikolai A. Vavilov, DENG Jiansong, Tatsuro Ito, Ian M. Wanless, Simone
Rinaldi and Jack H. Koolen  provided  some important  suggestions.
The  authors are  grateful to them. In 2012, when the  authors  read
an old paper on Latin  squares  and LYaPAS related to this article
written in Russian by Prof. Galina Borisovna Belyavskaya,  a  leading
researcher of the Institute of Mathematics and Computer Science of
the Academy of Sciences of the Republic of Moldova. Prof. Belyavskaya
was very kind to help and gave some very useful suggestions. Unfortunately,
Prof. Belyavskaya passed away in the morning of May 7, 2015, just
18 days after her 75'th birthday. We will  remember  this renowned
expert.

\bibliographystyle{amsplain}
\phantomsection\addcontentsline{toc}{section}{\refname}\bibliography{Ref}

\providecommand{\bysame}{\leavevmode\hbox to3em{\hrulefill}\thinspace}
\providecommand{\MR}{\relax\ifhmode\unskip\space\fi MR } 
\providecommand{\MRhref}[2]{%
  \href{http://www.ams.org/mathscinet-getitem?mr=#1}{#2}
} \providecommand{\href}[2]{#2} 
\end{document}